# Prototiles and Tilings from Voronoi and Delone cells of the Root Lattice $A_n$

Nazife Ozdes Koca[a)], Abeer Al-Siyabi[b)]

Department of Physics, College of Science, Sultan Qaboos University
P.O. Box 36, Al-Khoud, 123 Muscat, Sultanate of Oman,

Mehmet Koca[c)]

Department of Physics, Cukurova University, Adana, Turkey,

and

Ramazan Koc[d)]

Department of Physics, Gaziantep University, 27310, Gaziantep, Turkey

## ABSTRACT

We exploit the fact that two-dimensional facets of the Voronoi and Delone cells of the root lattice $A_n$ are identical rhombuses and equilateral triangles respectively. Their orthogonal projections onto the Coxeter plane display various rhombic and triangular prototiles including thick and thin rhombi of Penrose, Amman-Beenker tiles, Robinson triangles and Danzer triangles to name a few. We point out that the dihedral subgroup of order $2h$ involving the Coxeter element of order $h = n + 1$ of the Coxeter-Weyl group $a_n$ plays a crucial role for $h$-fold symmetric tilings of the Coxeter plane. After setting the general scheme we give samples of patches with 4, 5, 6, 7, 8 and 12-fold symmetries. The face centered cubic (f.c.c.) lattice described by the root lattice $A_3$ whose Wigner-Seitz cell is the rhombic dodecahedron projects, as expected, onto a square lattice with an $h = 4$ fold symmetry.

[a)]electronic-mail: nazife@squ.edu.om
[b)]electronic-mail: s21168@student.squ.edu.om
[c)]electronic-mail: mehmetkocaphysics@gmail.com ; retired professor
[d)]electronic-mail: koc@gantep.edu.tr

## 1. Introduction

Discovery of a 5-fold symmetric material (Shechtman et al., 1984) has led to growing interest in quasicrystallography. For a review see for instance (DiVincenzo & Steinhardt, 1991; Janot, 1993; Senechal, 2009). Aperiodic tilings of the plane with dihedral point symmetries or icosahedral symmetry in three dimensions have been the central research area of mathematicians and mathematical physicists to explain the quasicrystallography. For an excellent review see for instance (Baake & Grimm, 2013) and (Grünbaum & Shephard, 1987). There have been three major approaches for aperiodic tilings. The first class, perhaps, is the intuitive approach like Penrose tilings (Penrose, 1974 & 1978) which also exhibits the inflation-deflation technique developed later. The second approach is the projection technique of higher dimensional lattices onto lower dimensions pioneered by de Bruijn (de Bruijn, 1981) by projecting the 5-dimensional cubic lattice onto a plane orthogonal to one of the space diagonals of the cube. The cut-and-project technique works as follows. The lattice is partitioned into two complementing subspaces called $E_{\parallel}$ and $E_{\perp}$. Imagine a cylinder based on the component of the Voronoi cell in the subspace $E_{\perp}$. Project the lattice points lying in the cylinder into the space $E_{\parallel}$. Several authors (Duneau & Katz, 1985; Baake, Joseph, Kramer & Schlottmann, 1990; Chen, et al., 1998) have studied similar techniques. In particular, Baake et al. (Baake et al., 1990) exemplified the projection of the $A_4$ lattice. The group theoretical treatments of the projection of *n*-dimensional cubic lattices have been worked out by the references (Whittaker & Whittaker, 1987; Koca et al., 2015) and a general projection technique on the basis of dihedral subgroup of the root lattices has been proposed by Boyle & Steinhardt (Boyle & Steinhardt, 2016). In a recent article (Koca et al, 2018a) we pointed out that the 2-dimensional facets of the $A_4$Voronoi cell projects onto *thick* and *thin* rhombuses of the Penrose tilings. The third method is the model set technique initiated by Meyer (Meyer, 1972) and followed by Lagarias (Lagarias, 1996) and developed by Moody (Moody, 1997). For a detailed treatment see the reference (Baake & Grimm, 2013).

The Voronoi (Voronoi, 1908, 1909) and Delaunay (Delaunay, 1929, 1938a, b) cells of the root and weight lattices have been extensively studied in the inspiring book by Conway and Sloane (Conway & Sloane, 1988; chapter 21) and especially in the reference (Conway and Sloane, 1991). The reference (Deza & Grishukhin, 2004) contains detailed discussions and information about the Delone and Voronoi polytopes of the root lattices. Higher dimensional lattices have been worked out in (Engel, 1986) and the lattices of the root systems whose point groups are the Coxeter-Weyl groups have been studied extensively in the reference (Engel et al., 1994). Numbers of facets of the Voronoi and Delone cells of the root lattice have been also determined by a technique of decorated Coxeter-Dynkin diagrams (Moody & Patera, 1992). In a recent article (Koca et al., 2018b) we have worked out the detailed structures of the facets of the Voronoi and the Delone cells of the root and weight lattices of the $A_n$ and $D_n$ series. Basic information about the regular polytopes as the orbits of the Coxeter-Weyl groups have been worked out in the references (Coxeter, 1973) and (Grünbaum, 1967).

In this paper we point out that the facets of the Voronoi polytope of the root lattice $A_n$ are obtained from the Voronoi polytope of the cubic lattice $B_{n+1}$ by projection into n-dimensional Euclidean space along one of its space diagonals. It turns out that the facets of the Voronoi polytope of $A_n$ are rhombohedra in various dimensions such as rhombi in 2D, rhombohedra in 3D and the generalizations to higher dimensions. Projections of the 2D facets of the Voronoi polytope of $A_n$ lead to (n+1)-symmetric tilings of the Coxeter plane by a number of rhombi with different interior angles.

Similarly, we show that the 2D facets of the Delone cells of a given lattice $A_n$ are equilateral triangles; when projected into the Coxeter plane they lead to various triangles which tile the plane in an aperiodic manner. They include, among many new prototiles and patches, some well-known prototiles such as rhombus tilings by Penrose and Amman-Beenker and triangle tilings by Robinson and Danzer. The paper displays the lists of the rhombic and triangular prototiles in Table 1 and Table 2 originating from the projections of the Voronoi and Delone cells of the lattice $A_n$ onto the Coxeter plane respectively and illustrates some patches of the aperiodic tilings. Most of the prototiles and the tilings are displayed for the first time in what follows.

We organize the paper as follows. In Sec. 2 we introduce the basics of the root lattice $A_n$ via its Coxeter-Dynkin diagram, its Coxeter-Weyl group and projection of its Voronoi cell from that of $B_{n+1}$. The projections of the 2-dimensional facets of the Voronoi and Delone cells onto the Coxeter plane are studied in Section 3. Section 4 deals with the examples of periodic and aperiodic patches of tilings. Sec. 5 includes the concluding remarks. In appendix A we prove that the facets of the Voronoi cell of the 5D cubic lattice projects to those of the Voronoi cell of $A_4$.

## 2. The Root Lattice $A_n$ and its Coxeter-Weyl Group

We follow the notations of the reference (Conway and Sloane, 1991) which can be compared with the standard notations (Humphreys, 1990) and introduce some new ones when they are needed. The Coxeter-Dynkin diagram of $a_n$ describing the Coxeter-Weyl group and its extended diagram representing the affine Coxeter-Weyl group are shown in Fig.1.

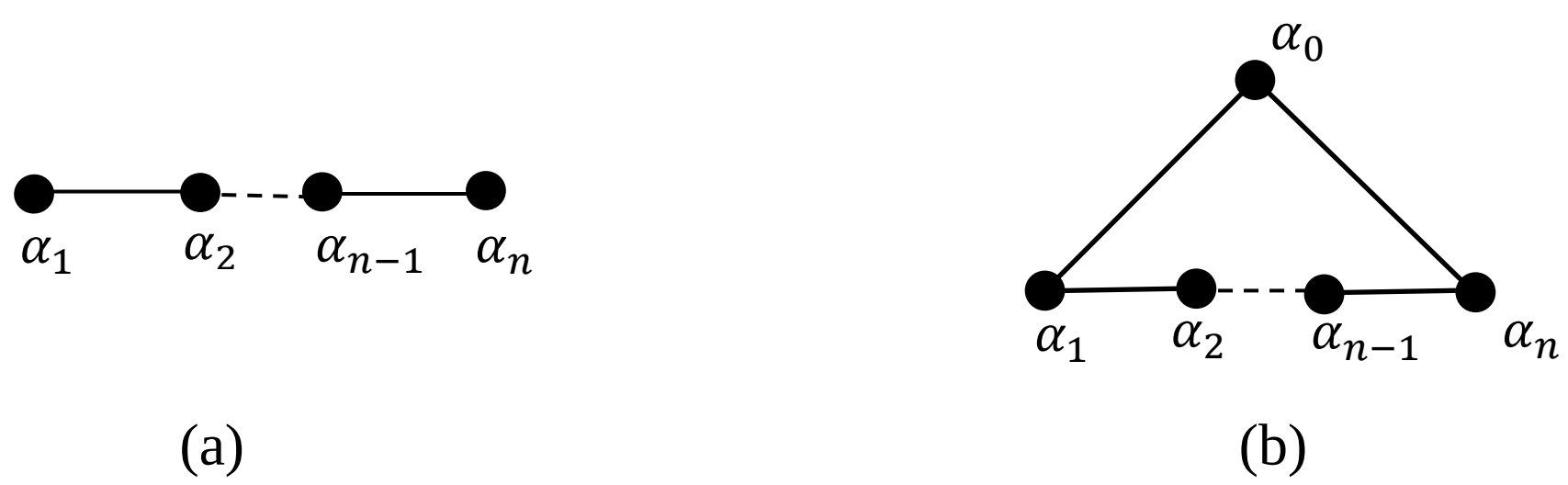

**Figure 1**
(a) Coxeter-Dynkin diagram of $a_n$, (b) Extended Coxeter-Dynkin diagram of $a_n$.

The nodes represent the lattice generating vectors (root vectors in Lie algebra terminology) $\alpha_i$ $(i = 1, 2, \dots, n)$ forming the Cartan matrix (Gram matrix of the lattice) by the relation

$$C_{ij} = \frac{2(\alpha_i, \alpha_j)}{(\alpha_j, \alpha_j)}\,. \qquad (1)$$

The fundamental weight vectors $\omega_i$ which generate the dual lattice $A_n^*$ satisfying the relation $(\alpha_i, \omega_j) = \delta_{ij}$ are given by the relations,

$$\omega_i = \sum_j (C^{-1})_{ij}\alpha_j, \quad \alpha_i = \sum_j C_{ij}\omega_j. \qquad (2)$$

The lattice $A_n$ is defined as the set of vectors $p = \sum_{i=1}^{n} b_i\,\alpha_i,\; b_i \in \mathbb{Z}$ and the weight lattice $A_n^*$ consists of the vectors $q = \sum_{i=1}^{n} c_i\,\omega_i,\; c_i \in \mathbb{Z}$ . Note that the root lattice is a sublattice of the weight lattice, $A_n \subset A_n^*$. Let $r_i, (i = 1, 2, \dots, n)$ denotes the reflection generator with respect to the hyperplane orthogonal to the simple root $\alpha_i$ which operates on an arbitrary vector $\lambda$ as

$$r_i\lambda = \lambda - \frac{2(\lambda,\alpha_i)}{(\alpha_i,\alpha_i)}\alpha_i \,. \tag{3}$$

The reflection generators generate the Coxeter-Weyl group $< r_1, r_2, \dots, r_n \,|(r_i\, r_j)^{m_{ij}} = 1 >$. Adding another generator, usually denoted by $r_0$, describing the reflection with respect to the hyperplane bisecting the highest weight vector $(\omega_1 + \omega_n)$ we obtain the *Affine Coxeter group*, the infinite discrete group denoted by $< r_0, r_1, r_2, \dots, r_n >$. The relations between the Voronoi and the Delone cells of the lattices $A_n$ and $A_n^*$ have been studied by L. Michel (Michel, 1995, 2001). An arbitrary group element of the Coxeter-Weyl group will be denoted by $W(a_n)$ and the orbit of an arbitrary vector $q$ will be defined as $W(a_n)q =: (\sum_{i=1}^{n} c_i\, \omega_i)_{a_n}$, $c_i \in \mathbb{Z}$. With this notation the root polytope will be denoted either by $(10 \dots 01)_{a_n}$ or simply by $(\omega_1 + \omega_n)_{a_n}$. The dual polytope of the root polytope (Koca et al., 2018b) is the Voronoi cell $V(0)$ which is the union of the orbits of the fundamental polytopes $(\omega_i)_{a_n}, (i = 1, 2, \dots, n)$,

$$(\omega_1)_{a_n} \cup (\omega_2)_{a_n} \dots \cup (\omega_n)_{a_n}. \tag{4}$$

Each fundamental polytope $(\omega_i)_{a_n}$ of $A_n$ is a copy of one of the Delone polytopes. Delone polytopes tile the root lattice in such a way that each polytope centralizes one vertex of the Voronoi cell. For example, the set of vertices $(\omega_1)_{a_n} + (\omega_n)_{a_n}$ represent the $2(n+1)$ simplexes centered around the vertices $(\omega_1)_{a_n}$ and $(\omega_n)_{a_n}$ of the Voronoi cell $V(0)$. Similarly, $(\omega_2)_{a_n} + (\omega_{n-1})_{a_n}$constitutes the vertices of the ambo-simplexes (see the reference (Conway and Sloane, 1991) for the definition) centralizing the vertices $(\omega_2)_{a_n}$ and the $(\omega_{n-1})_{a_n}$ of the Voronoi cell $V(0)$ and so on.

It is a general practice to work with the set of orthonormal vectors $l_i, (i = 1, 2, \dots, n+1), (l_i, l_j) = \delta_{ij}$ to define the root and weight vectors. The root vectors can be defined as $\alpha_i = l_i - l_{i+1}, (i = 1, 2, \dots, n)$ which are permuted by the generators as $r_i\colon l_i \leftrightarrow l_{i+1}$ implying that the Coxeter-Weyl group $S_{n+1}$ is isomorphic to the symmetric group of order $(n+1)!$ permuting the $n+1$ vectors. The fundamental weights then are given as

$$\omega_1 = \frac{1}{n+1}(n\, l_1 - l_2 - l_3 - \dots - l_{n+1});$$

$$\omega_2 = \frac{1}{n+1}\big((n-1)\, l_1 + (n-1)\, l_2 - 2\, l_3 - \dots - 2\, l_{n+1}\big);$$

$$\vdots$$

$$\omega_n = \frac{1}{n+1}(\, l_1 + l_2 + \dots + l_n - n\, l_{n+1}). \tag{5}$$

It is even more convenient to introduce the non-orthogonal set of vectors

$$k_i =: -\frac{l_0}{n+1} + l_i, \quad l_0 =: l_1 + l_2 + \dots + l_{n+1}, \qquad k_1 + k_2 + \dots + k_{n+1} = 0. \tag{6}$$

The vectors $k_i$ are more useful since they represent the vertices of the $n$-simplex $(\omega_1)_{a_n}$ and the edges of the Voronoi cell are parallel to the vectors $k_i$ . In terms of the fundamental weights they are given by

$$k_1 = \omega_1, k_2 = \omega_2 - \omega_1, k_3 = \omega_3 - \omega_2, \ldots, k_n = \omega_n - \omega_{n-1},\ k_{n+1} = -\omega_n. \quad (7)$$

Their magnitudes are all the same $(k_i, k_i) = \frac{n}{n+1}$ and the angle between any pair is $\theta = \cos^{-1}\left(-\frac{1}{n}\right)$. They are orthogonal to the vector $l_0$, $(k_i,\ l_0) = 0$ so that they define the $n-$dimensional hyperplane orthogonal to the vector $l_0$. Following (5) and (6) vertices of the orbits of the fundamental weights

$$(\omega_i)_{a_n} = (k_1 + k_2 + \cdots + k_i)_{a_n},\ (i = 1, 2, \ldots, n) \quad (8)$$

are obtained as the permutations of the vectors in (8) and each orbit of the fundamental weight represents a copy of the Delone polytope. It is clear that the generators of the Coxeter-Weyl group permute the vectors $k_i$ as $r_i$: $k_i \leftrightarrow k_{i+1}$. The number of vertices of the Voronoi cell of $A_n$ is then given by

$$\sum_{i=1}^{n} \binom{n+1}{i} = 2^{n+1} - 2. \quad (9)$$

The Voronoi cell $V(0)$ of $A_n$ can also be obtained as the projection along the space diagonal vector $\pm\, l_0$ of the Voronoi cell of the cubic lattice in $n+1$ dimensions with vertices $\frac{1}{2}(\pm\, l_1 \pm\ l_2 \pm \cdots \pm\ l_{n+1})$.
Substituting $l_i =: \frac{l_0}{n+1} + k_i$ it is straightforward to show that the vertices of the Voronoi cell of $B_{n+1}$ project into the vertices of the Voronoi cell of $A_n$ represented by the union of the fundamental polytopes in n-dimensional space given in (4).
This follows from the projections

$$\pm\frac{1}{2} l_0 \to 0,$$

$$\pm\frac{1}{2}(-l_1 + l_2 + \cdots + l_{n+1}) \to \mp\, k_1 = \mp\omega_1,$$

$$\pm\frac{1}{2}(-l_1 - l_2 + \cdots + l_{n+1}) \to \mp\,(k_1 + k_2) = \mp\omega_2, \quad (10)$$

$$\vdots$$

$$\pm\frac{1}{2}(-l_1 - l_2 - \cdots - l_n + l_{n+1}) \to \mp\,(k_1 + k_2 + \cdots + k_n) = \pm\,(k_{n+1}) = \mp\omega_n.$$

Let us consider the root system $l_i - l_j = k_i - k_j$ ; the hyperplanes bisecting the segment between the origin and the root system are represented by the vectors $\frac{1}{2}(k_i - k_j)$. The intersections of these hyperplanes are the vertices of the Voronoi cell $V(0)$ of the lattice $A_n$ which are the union of the orbits of the highest weight vectors given in (4). Actually the vectors $\frac{1}{2}(k_i - k_j)$ are the centers of the $n-1$ dimensional rhombohedral facets of the Voronoi cell $V(0)$. The $2D$ facets (squares), $3D$ (cubes) and higher dimensional cubic facets of the Voronoi cell of $B_{n+1}$ project into the rhombi, rhombohedra and higher dimensional rhombohedra representing the facets of the Voronoi cell of the $V(0)$ of $A_n$. Appendix A discusses further aspects of this projection for two special cases $B_3 \to A_2$ and $B_5 \to A_4$.

Two dimensional faces of the Delone polytopes are equilateral triangles. To give an example let us take the vertex $\omega_1 = k_1$. The set of vertices $(\omega_1)_{a_n} = \{k_i\}, i = 1, 2, \dots, n+1$ Constitute the vertices of the *n*-simplex. One can generate an equilateral triangle by applying the subgroup $< r_1, r_2 >$ on $\omega_1$ where the edges are given by the root vectors $k_1 - k_2$, $k_2 - k_3$, $k_3 - k_1$. Other faces and the corresponding edges are generated by applying the group elements $W(a_n)$ leading to the set of edges $k_i - k_j$, a vector in the root system. The argument is valid for any other Delone polytope $(\omega_i)_{a_n}$. There is always a group element relating the vector $k_i$ to any other vertex $k_j$.To give a further example let us consider the lattice $A_3$. The Delone polytopes $(\omega_1)_{a_3} = -(\omega_3)_{a_3}$ represent two tetrahedra where the vertices belong to the set $\{k_i\}$ or $\{-k_i\}$ . The six edges of the four triangles of each tetrahedra are represented by the vectors $k_i - k_j, i \neq j = 1,2,3,4$ belonging to the root system. Similarly, the orbit $(\omega_2)_{a_3}$ is an octahedron with vertices $k_i + k_j, i \neq j = 1,2,3,4$ where the edges of the eight triangular faces are represented by the vectors $k_i - k_j, i \neq j$ the elements of the root system. For a general proof one can take the vertex in (8) and apply the Coxeter $-$ Weyl group $a_n$ to generate the nearest vertices forming triangles associated with the vertex (8). These arguments imply that the set of edges of an arbitrary two dimensional face of a Delone polytope can be represented by the root vectors $k_i - k_j$, $k_j - k_l$, $k_l - k_i$, $i \neq j \neq l = 1, 2, \dots, n+1$.

## 3. Projections of the faces of the Voronoi and Delone cells

In this section we define the Coxeter plane and determine the dihedral subgroup of the Coxeter-Weyl group acting in this plane. The Coxeter plane is generally defined by the eigenvector $v = (\zeta, \zeta^2, \dots, \zeta^n, 1)$ corresponding to the eigenvalue $\zeta = e^{\frac{i2\pi}{h}}$ of the Coxeter element $R = r_1 r_2 r_3 \dots r_n$ of the group $a_n$ which permutes the vectors $k_1 \to k_2 \to \cdots k_{n+1} \to k_1$ in the cyclic order. Here $h = n+1$ is the Coxeter number. Going to the real space, the Coxeter element acts on the plane $E_\parallel$ spanned by the unit vectors $u_1 = \frac{1}{\sqrt{2h}}(v + \bar{v})$ and $u_2 = \frac{i}{\sqrt{2h}}(v - \bar{v})$ . The components of the vectors $k_j$ in the Coxeter plane is then represented by the complex number $(k_j)_\parallel = ce^{\frac{ij2\pi}{h}}$, where $j = 1,2, \dots, n+1$ and $c = \sqrt{\frac{2}{h}}$.

In (Koca, et al., 2014, 2015) we have introduced an equivalent definition of the Coxeter plane through the eigenvalues and eigenvectors of the Cartan matrix of the root system of the lattice $A_n$. The eigenvalues and eigenvectors of the Cartan matrix of the group $a_n$ can be written as

$$\lambda_i = 2(1 + \cos\frac{m_i\pi}{h}); \quad X_{ji} = (-1)^{j+1} \sin\frac{j\, m_i\pi}{h}, \quad i, j = 1, 2, \dots, n, \qquad (11)$$

where $m_i = 1, 2, \dots, n$ are the Coxeter exponents. We define the orthonormal set of vectors subject to an arbitrary orthogonal transformation by

$$\hat{x}_i = \frac{1}{\sqrt{\lambda_i}} \sum_{j=1}^{n} \alpha_j X_{ji}, \ \left(\hat{x}_i, \hat{x}_j\right) = \delta_{ij}, \ i = 1, 2, \dots, n. \tag{12}$$

An equivalent definition of the Coxeter plane $E_{\parallel}$ is then given by the pair of vectors $\hat{x}_1$ and $\hat{x}_n$ and the rest of the vectors in (12) determine the orthogonal space $E_{\perp}$. With a suitable choice of the simple roots $\alpha_j$ the unit vectors $\hat{x}_1$ and $\hat{x}_n$ can be associated with the unit vectors $u_1$ and $u_2$ respectively.

The Coxeter element $R = r_1 r_2 r_3 \dots r_n$ can be written as the product of two reflection generators $< R_1, R_2 \,|(R_1 R_2)^h = 1 >$ defining a dihedral subgroup of order $2h$ of the group $a_n$ (Carter, 1972; Humphreys, 1990). Therefore, the residual symmetry after projection is the dihedral group acting on the Coxeter plane.

We will drop the overall factor $c$ as it has no significant meaning in our further discussions. It is then clear that all possible scalar products $((k_j)_{\parallel}, (k_l)_{\parallel}) = \cos\left(\frac{(j-l)2\pi}{h}\right)$ will determine the nature of the projected rhombi in the Coxeter plane $E_{\parallel}$. The tiles projected from the Delone cells will be triangles whose edge lengths are the magnitudes of the vectors $(k_i - k_j)_{\parallel}, (k_j - k_l)_{\parallel}, (k_l - k_i)_{\parallel}$ and, after deleting the common factor $c$, they will read respectively

$$2\left|\sin\left(\frac{(i-j)\pi}{h}\right)\right|, 2\left|\sin\left(\frac{(j-l)\pi}{h}\right)\right|, 2\left|\sin\left(\frac{(l-i)\pi}{h}\right)\right|, \tag{13}$$

exhausting all possible projections. For each set of integral numbers $i \neq j \neq l = 1, 2, \dots, n+1$ we may assume without loss of generality that the numbers can be ordered as $i > j > l$ and the corresponding angles of the edges in (13) can be represented as $\frac{(i-j)\pi}{h}, \frac{(j-l)\pi}{h}, \pi - \frac{(i-l)\pi}{h}$. Defining these angles by $\left(\frac{n_1\pi}{h}, \frac{n_2\pi}{h}, \frac{n_3\pi}{h}\right)$, $n_i \in \mathbb{N}$ we obtain $n_1 + n_2 + n_3 = h$. Some edges from the Voronoi cell or Delone cell overlap after projection forming a degenerate rhombi or triangle respectively which has no effect on tilings of the Coxeter plane.

In what follows we will illustrate these results with some examples. The rhombic and triangular prototiles are classified in Table 1 and Table 2 respectively.

**Table 1**
Classification of rhombic prototiles projected from Voronoi cells of the root lattice $A_n$ (rhombuses with angles $\left(\frac{2\pi m}{h}, \pi - \frac{2\pi m}{h}\right), m \in N$).

| Root Lattice | $h$ | # of prototiles | Rhombuses with pairs of angles |
|---|---|---|---|
| $A_3$ | 4 | 1 | $\left(\frac{\pi}{2}, \frac{\pi}{2}\right)$, (square) |
| $A_4$ | 5 | 2 | $\left(\frac{2\pi}{5}, \frac{3\pi}{5}\right), \left(\frac{4\pi}{5}, \frac{\pi}{5}\right)$, (Penrose's thick and thin rhombuses) |
| $A_5$ | 6 | 1 | $\left(\frac{\pi}{3}, \frac{2\pi}{3}\right)$ |
| $A_6$ | 7 | 3 | $\left(\frac{2\pi}{7},\frac{5\pi}{7}\right),\left(\frac{4\pi}{7},\frac{3\pi}{7}\right),\left(\frac{6\pi}{7},\frac{\pi}{7}\right)$ |
| $A_7$ | 8 | 2 | $\left(\frac{\pi}{4},\frac{3\pi}{4}\right),\left(\frac{\pi}{2},\frac{\pi}{2}\right)$, (Amman-Beenker tiles) |
| $A_8$ | 9 | 4 | $\left(\frac{2\pi}{9},\frac{7\pi}{9}\right),\left(\frac{4\pi}{9},\frac{5\pi}{9}\right),\left(\frac{6\pi}{9},\frac{3\pi}{9}\right),\left(\frac{8\pi}{9},\frac{\pi}{9}\right)$ |
| $A_9$ | 10 | 2 | $\left(\frac{2\pi}{5},\frac{3\pi}{5}\right),\left(\frac{\pi}{5},\frac{4\pi}{5}\right)$, (Penrose' thick and thin rhombuses) |
| $A_{10}$ | 11 | 3 | $\left(\frac{2\pi}{11},\frac{9\pi}{11}\right),\left(\frac{4\pi}{11},\frac{7\pi}{11}\right),\left(\frac{6\pi}{11},\frac{5\pi}{11}\right)$ |
| $A_{11}$ | 12 | 3 | $\left(\frac{\pi}{6},\frac{5\pi}{6}\right),\left(\frac{\pi}{3},\frac{2\pi}{3}\right),\left(\frac{\pi}{2},\frac{\pi}{2}\right)$ |

**Table 2**
Classification of triangular prototiles with angles $\left(\frac{n_1\pi}{h},\frac{n_2\pi}{h},\frac{n_3\pi}{h}\right)$ projected from Delone cells of the root lattice $A_n$.

| Root Lattice | $h$ | # of prototiles | Triangles denoted by triple natural numbers $(n_1,n_2,n_3)$; $n_1+n_2+n_3=h$ |
|---|---|---|---|
| $A_2$ | 3 | 1 | $(1,1,1)$ |
| $A_3$ | 4 | 1 | $(1,1,2)$, (right triangle) |
| $A_4$ | 5 | 2 | $(1,1,3),(1,2,2)$, (Robinson triangles) |
| $A_5$ | 6 | 3 | $(1,1,4),(1,2,3),(2,2,2)$ |
| $A_6$ | 7 | 4 | $(1,1,5),(1,2,4),(1,3,3),(2,2,3)$, (Danzer triangles) |
| $A_7$ | 8 | 5 | $(1,1,6),(1,2,5),(1,3,4),(2,2,4),(2,3,3)$ |
| $A_8$ | 9 | 7 | $(1,1,7),(1,2,6),(1,3,5),(1,4,4),(2,2,5),(2,3,4),(3,3,3)$ |
| $A_9$ | 10 | 8 | $(1,1,8),(1,2,7),(1,3,6),(1,4,5),(2,2,6),(2,3,5),(2,4,4),(3,3,4)$ |
| $A_{10}$ | 11 | 10 | $(1,1,9),(1,2,8),(1,3,7),(1,4,6),(1,5,5),(2,2,7),(2,3,6),(2,4,5),(3,3,5),(3,4,4)$ |
| $A_{11}$ | 12 | 12 | $(1,1,10),(1,2,9),(1,3,8),(1,4,7),(1,5,6),(2,2,8),(2,3,7),(2,4,6),(2,5,5),(3,3,6),(3,4,5),(4,4,4)$ |

## 4. Examples of prototiles and patches of tilings

In what follows we discuss a few $h$-fold symmetric tilings with rhombuses and triangles.

### 4.1a. Projection of the Voronoi cell of $A_3$

The Voronoi cell of the root lattice $A_3$ , commonly known as the f.c.c. lattice, is the rhombic dodecahedron (the Wigner-Seitz cell) which tessellates 3-dimensional Euclidean space. The vertices of the rhombic dodecahedron can be represented as the union of two tetrahedra and one octahedron:

$$(\omega_1)_{a_3} \cup (\omega_2)_{a_3} \cup (\omega_3)_{a_3} = \{(\pm k_1, \pm k_2, \pm k_3, \pm k_4);$$

$$(k_1 + k_2,\ k_2 + k_3,\ k_3 + k_4,\ k_4 + k_1,\ k_1 + k_3,\ k_2 + k_4)\}. \qquad (14)$$

A typical rhombic face of the rhombic dodecahedron is determined by the set of four vertices $k_1 = \omega_1,\ k_1 + k_2 = \omega_2,\ k_1 + k_3 = r_2\omega_2,\ \omega_3 = k_1 + k_2 + k_3 = -k_4$.

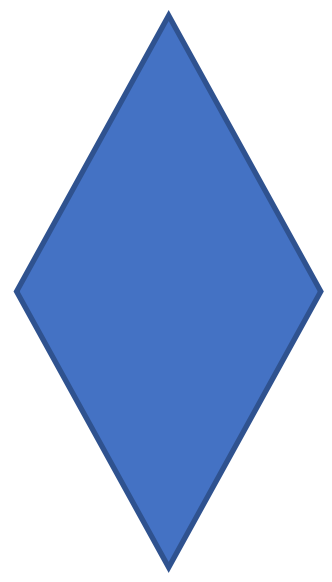

**Figure 2**
A typical rhombic face of the rhombic dodecahedron (note that edges are represented by the vectors $k_2$ and $k_3$ with an angle $109.5^0$).

The first four vectors in (14), $(\omega_1)_{a_3} = \{k_1, k_2, k_3, k_4\}$ represent the vertices of the first tetrahedron and the permutation group $S_4 \cong W(a_3)$ is isomorphic to the tetrahedral group of order 24. The set of vectors $(-k_1, -k_2, -k_3, -k_4)$ represent the second tetrahedron $(\omega_3)_{a_3}$. Actually, the Dynkin diagram symmetry $(k_1 \leftrightarrow -k_4)$ extends the tetrahedral group to the octahedral group implying that two tetrahedra form a cube. Since the Voronoi cell is dual to the root polytope it is face transitive and invariant under the octahedral group of order 48. The last six vectors in (14) represent the vertices of an octahedron which split as $4 + 2$ under the dihedral subgroup $D_4$ of order 8. The last two vectors $k_1 + k_3$ and $k_2 + k_4 = -(k_1 + k_3)$ form a doublet under the dihedral subgroup and they are perpendicular to the Coxeter plane $E_\parallel$ determined by the pair of vectors $(\hat{x}_1, \hat{x}_3)$. The projected components of the vectors $k_i$ in the plane $E_\parallel$ are given by

$$(k_1)_\parallel = (0, 1),\ (k_2)_\parallel = (-1, 0),\ (k_3)_\parallel = (0, -1),\ (k_4)_\parallel = (1, 0). \qquad (15)$$

Therefore, projection of the rhombus in Fig. 2 is a square of unit length determined by the vectors $(k_2)_{\|} = (-1, 0)$ and $(k_3)_{\|} = (0, -1)$. Projection of the Voronoi cell onto the plane $E_{\|}$ is shown in Fig. 3 which is invariant under the dihedral group $D_4$.We can pick up any two vectors in (15) say $(k_1)_{\|}$ and $(k_2)_{\|}$ then the lattice in $E_{\|}$ is a square lattice with a general vector $q_{\|} = m_1(k_1)_{\|} + m_2(k_2)_{\|}$ with $m_1, m_2 \in \mathbb{Z}$ .

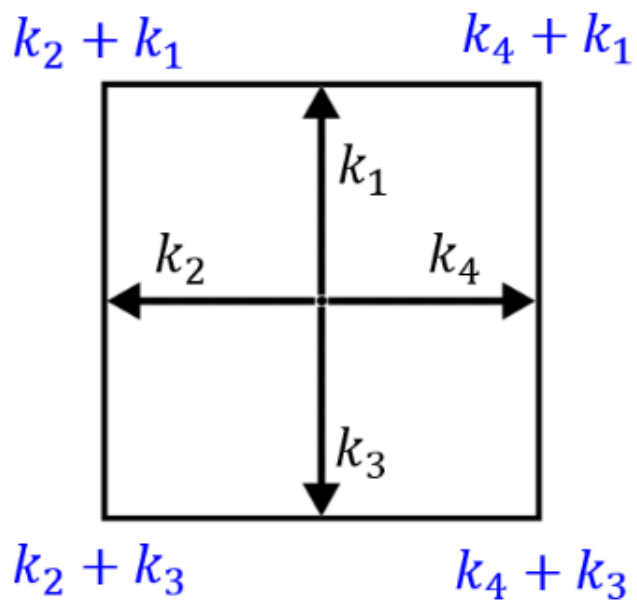

**Figure 3**
Projection of the Wigner-Seitz cell (rhombic dodecahedron) onto the Coxeter plane $E_{\|}$.

**4.1b. Projection of the Delone cells of the root lattice $A_3$**

The lattice vector of the lattice $A_3$ can be written as $p = \sum_{i=1}^{3} b_i \alpha_i = \sum_{i=1}^{3} n_i k_i$ with $b_i, n_i \in \mathbb{Z}$ such that $\sum_{i=1}^{3} n_i$ =even. This shows that the root lattice is a sublattice of the weight lattice. Equilateral triangles of the Delone cells project as *right triangles*, any two constitute a square of length $\sqrt{2}$ . One can check that the projected Delone cells form a square lattice with a general vector $p_{\|} = m_1(k_1)_{\|} + m_2(k_2)_{\|}$, with $m_1 + m_2$ =even integer as shown in Fig. 4.

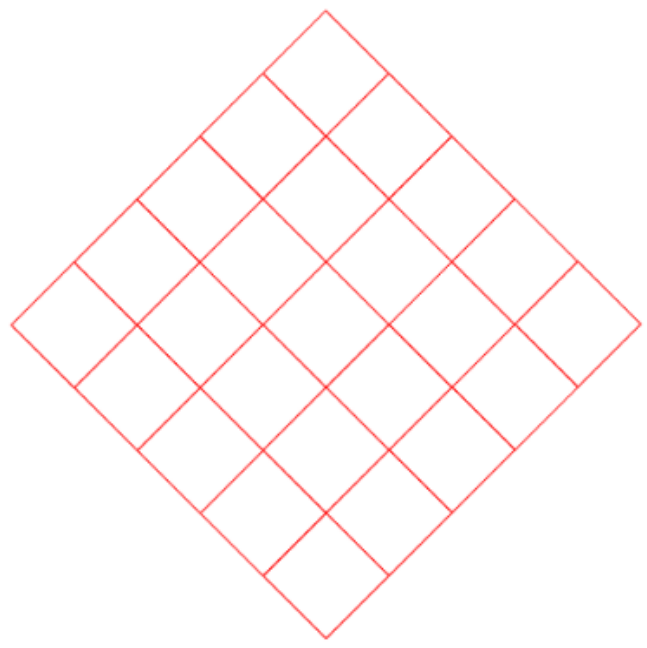

**Figure 4**
Projection of the root lattice as a square lattice made of right triangles (not shown).

It is clear that it is a sublattice of the square lattice obtained from the weight lattice and it is rotated by $45^0$ with respect to the projected square lattice of the Voronoi cell. With this example we have obtained a sublattice whose unit cell is invariant under the dihedral group $D_4$. This is expected because the dihedral group is a crystallographic group. Two overlapping square lattices have been depicted in Fig. 5.

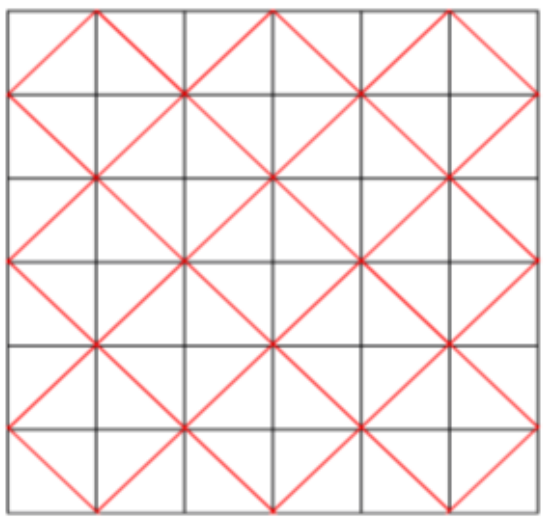

**Figure 5**
Projections of the lattices $A_3{}^*$ and $A_3$ leading to two overlapping square lattices, $A_3$ is the sub lattice of the other.

### 4.2a. Projection of the Voronoi cell of the root lattice $A_4$

We would like to discuss this projection leading to the Penrose tiling which can be obtained from the cut-and-project technique or by the substitution technique directly. The Voronoi cell of the root lattice tiles the 4-dimensional Euclidean reciprocal space which is the weight lattice $A_4{}^*$ represented by a general vector $q = \sum_{i=1}^{4} c_i\, \omega_i = \sum_{i=1}^{5} m_i\, k_i = \sum_{i=1}^{4} n_i\, k_i$ with $c_i, m_i, n_i \in \mathbb{Z}$. Union of the Delone polytopes

$$(\omega_1)_{a_4} \cup (\omega_2)_{a_4} \cup (\omega_3)_{a_4} \cup (\omega_4)_{a_4} \quad (16)$$

constitutes the Voronoi cell $V(0)$ with $30 = 5 + 10 + 10 + 5$ vertices. It is a polytope comprised of rhombohedral facets. A typical 3-facet of the Voronoi cell is a rhombohedron with six rhombic faces (Koca et al., 2018b) as depicted in Fig. 6. Its vertices are given by

$$\left(k_1,\ k_1 + k_2,\ k_1 + k_3,\ k_1 + k_4, -k_5, -(k_2 + k_5), -(k_3 + k_5), -(k_4 + k_5)\right)$$

and its center is at $\frac{1}{2}(k_1 - k_5)$.

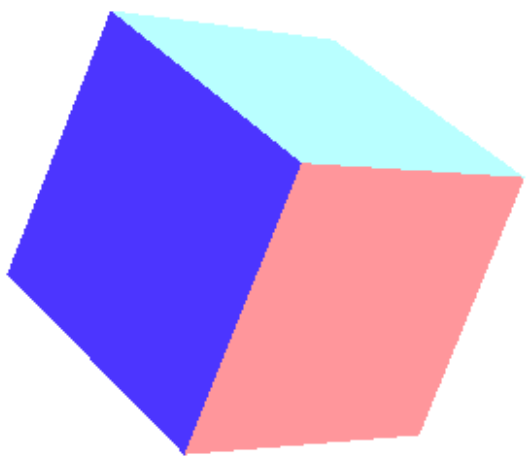

**Figure 6**
A 3-dimensional rhombohedron, the rhombohedral facet of the Voronoi cell of the root lattice $A_4$.

The Voronoi cell is formed by 20 such rhombohedra whose rhombic 2-faces are generated by the pairs of vectors satisfying the scalar product $(k_i,\ k_j) = -\frac{1}{5}, i \neq j = 1, \ldots, 5$ and $(k_i,\ k_i) = \frac{4}{5}$. It is perhaps useful to discuss some of the technicalities of the projection with this example. The components of the vectors in the plane $E_\parallel$ $(k_j)_\parallel = (l_j)_\parallel \Rightarrow \zeta^j = e^{\frac{ij\pi}{5}},\ j = 1, \ldots, 5$ satisfying $\sum_{i=1}^{5} \zeta^j = 0$ form a pentagon while the vectors $(k_j)_\perp$ form a pentagram since the angle between two successive vectors $(k_j)_\perp$ is $\frac{4\pi}{5}$.

A matrix representation of the Coxeter element $R$ in 4-dimensional space is given by

$$R = \begin{pmatrix} \cos\frac{2\pi}{5} & -\sin\frac{2\pi}{5} & 0 & 0 \\ \sin\frac{2\pi}{5} & \cos\frac{2\pi}{5} & 0 & 0 \\ 0 & 0 & \cos\frac{4\pi}{5} & -\sin\frac{4\pi}{5} \\ 0 & 0 & \sin\frac{4\pi}{5} & \cos\frac{4\pi}{5} \end{pmatrix}. \quad (17)$$

A general vector of the weight lattice $A_4{}^*$ projects as $q_\parallel = \sum_{i=1}^{5} n_i (k_i)_\parallel$ . Note that the rhombic 2-faces of the Voronoi cell project onto the plane $E_\parallel$ as two Penrose rhombuses, *thick* and *thin*, for the acute angles between any pair of the vectors $(k_j)_\parallel \Rightarrow e^{\frac{ij2\pi}{5}}$ are either $72^0$or $36^0$.

De Bruijn (de Bruijn, 1981) proved that the projection of the 5-dimensional cubic lattice onto $E_\parallel$ leads to the Penrose rhombus tiling of the plane with thick and thin rhombuses. The same technique is extended to an arbitrary cubic lattice (Whittaker & Whittaker, 1987). It is not a surprise that we obtain the same tiling from the $A_4$ lattice because Voronoi cell of $A_4$ is obtained as the projection of the Voronoi cell of the cubic lattice $B_5$ (Koca et al., 2015) (see also Appendix A for further details). The Coxeter-Weyl group $W(b_5)$ is of order $2^5.5!$ and its Coxeter number is $h = 10$. Normally one expects a 10-fold symmetric tiling from the projection of the 5-dimensional cubic lattice. The Voronoi cell of $B_5$ has 32 vertices $\frac{1}{2}(\pm l_1 \pm l_2 \pm l_3 \pm l_4 \pm l_5)$ which decompose under the Coxeter-Weyl group $W(a_4)$ as

$$32 = 1 + 1 + 5 + 10 + 10 + 5. \quad (18)$$

The first two terms in (18) are represented by the vectors $\pm\frac{1}{2} l_0 \rightarrow 0$ in 4-space and the remaining vertices of 5-dimensional cube decompose as $5 + 10 + 10 + 5$ representing the vertices of the Voronoi cell $(\omega_1)_{a_4} \cup (\omega_2)_{a_4} \cup (\omega_3)_{a_4} \cup (\omega_4)_{a_4}$. In terms of the vectors $l_i$, the vectors $k_i$ can also be written as (compare with (6))

$$\begin{aligned} k_1 - \frac{3}{10} l_0 &= \frac{1}{2}( l_1 - l_2 - l_3 - l_4 - l_5), \\ k_2 - \frac{3}{10} l_0 &= \frac{1}{2}(- l_1 + l_2 - l_3 - l_4 - l_5), \\ k_3 - \frac{3}{10} l_0 &= \frac{1}{2}(-l_1 - l_2 + l_3 - l_4 - l_5), \\ k_4 - \frac{3}{10} l_0 &= \frac{1}{2}(- l_1 - l_2 - l_3 + l_4 - l_5), \\ k_5 - \frac{3}{10} l_0 &= \frac{1}{2}(-l_1 - l_2 - l_3 - l_4 + l_5). \end{aligned} \quad (19)$$

This implies that the set of vertices of the 5-dimensional cube on the right of equation (19) represent a 4-simplex; any three $k_i$ represent the vertices of an equilateral triangle, any four constitute a tetrahedron. The 4-simplex is one of the Delone polytopes. Similarly, pairwise and triple combinations of the vectors $k_i$ constitute the Delone cells composed of octahedra and tetrahedra as 3-dimensional facets. De Bruijn proved that the integers $n_i$ in the expression $q_{\parallel} = \sum_{i=1}^{5} n_i (k_i)_{\parallel}$ take values $\sum_{i=1}^{5} n_i \in \{1, 2, 3, 4\}$. The vectors in the direction $(\pm k_i)_{\parallel}$ count positive and negative depending on its sign. A patch of Penrose tiling with the numbering of vertices is shown in Fig. 7 follows from the cut-and-project technique. Note that the projection of the Voronoi cell onto $E_{\parallel}$ form four intersecting pentagons defined by de Bruijn and denoted in our notation by

$$V_1 = \left((\omega_1)_{a_4}\right)_{\parallel} = -V_4 = -\left((\omega_4)_{a_4}\right)_{\parallel},\ V_2 = \left((\omega_2)_{a_4}\right)_{\parallel} = -V_3 = -\left((\omega_3)_{a_4}\right)_{\parallel}.$$

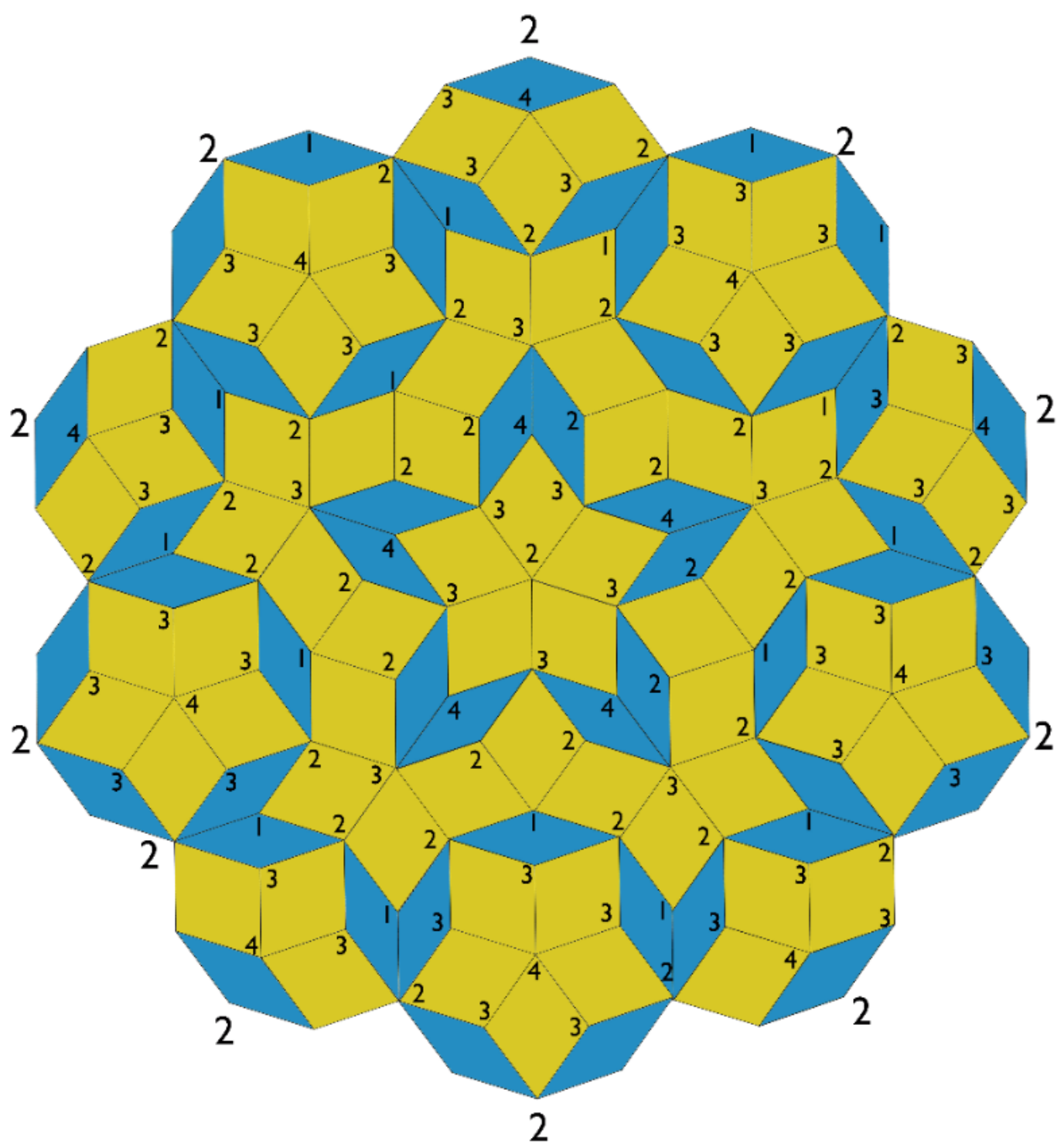

**Figure 7**
A patch of the Penrose rhombic tiling obtained by the cut-and-project method from the lattice $A_4{}^*$. The four types of vertices are distinguished by numbers as stated in the text.

### 4.2b. Projection of the root lattice $A_4$ by Delone cells

The Delone polytopes $(\omega_1)_{a_4} = -(\omega_4)_{a_4}, (\omega_2)_{a_4} = -(\omega_3)_{a_4}$, tile the root lattice such that the centers of the Delone cells correspond to the vertices of the Voronoi cells. Take, for example, the 4-simplex $(\omega_1)_{a_4} = \{k_1, k_2, k_3, k_4, k_5\}$ and $(\omega_4)_{a_4} = \{-k_1, -k_2, -k_3, -k_4, -k_5\}$. Adding vertices of these two simplexes one obtains 10 simplexes of Delone cells centered around 10 vertices of the Voronoi cell $V(0)$. If we add the vector $k_1$ to the vertices of the Delone cell $(\omega_4)_{a_4}$ we obtain the vertices of the Delone cell $\{0, k_1 - k_2, k_1 - k_3, k_1 - k_4, k_1 - k_5\}$ whose center is at the vertex $k_1$ of the Voronoi cell $V(0)$. When the vertices of the polytope $(\omega_2)_{a_4} = \{k_i + k_j, i \neq j = 1, 2, \dots, 5\}$ are added to the vertices of $(\omega_3)_{a_4}$ we obtain 20 Delone cells centered around the vertices of $(\omega_2)_{a_4}$ and $(\omega_3)_{a_4}$. For example when we add the vector $k_1 + k_2$ to the vertices of $(\omega_3)_{a_4}$we obtain 10 vertices of the Delone cell centered around the vertex $k_1 + k_2$ of the Voronoi cell $V(0)$ as

$$\{0,\ k_1 - k_3,\ k_1 - k_3 + k_2 - k_4,\ k_1 - k_4 + k_2 - k_5,\ k_2 - k_5,\ k_2 - k_3,\ k_1 - k_4,\ k_1 - k_3 + k_2 - k_5,\ k_2 - k_4,\ k_1 - k_5\}. \quad (20)$$

Note that the vertices of the Delone cells are either the root vectors or their linear combinations as expected.

Using (13) for $n = 4$ and $j \neq k \neq l = 1, 2, \ldots, 5$ we obtain two isosceles Robinson triangles with edges $(2\sin\frac{\pi}{5}, 2\sin\frac{\pi}{5}, 2\sin\frac{2\pi}{5})$ and $(2\sin\frac{2\pi}{5}, 2\sin\frac{2\pi}{5}, 2\sin\frac{4\pi}{5})$. The darts and kites obtained from the Robinson triangles are shown in Fig. 8.

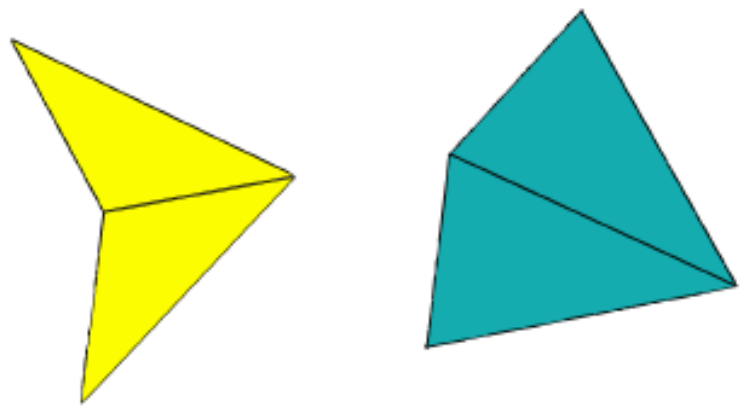

**Figure 8**
Darts and kites obtained from Robinson triangles.

A patch of Penrose tiling by darts and kites is shown Fig. 9.

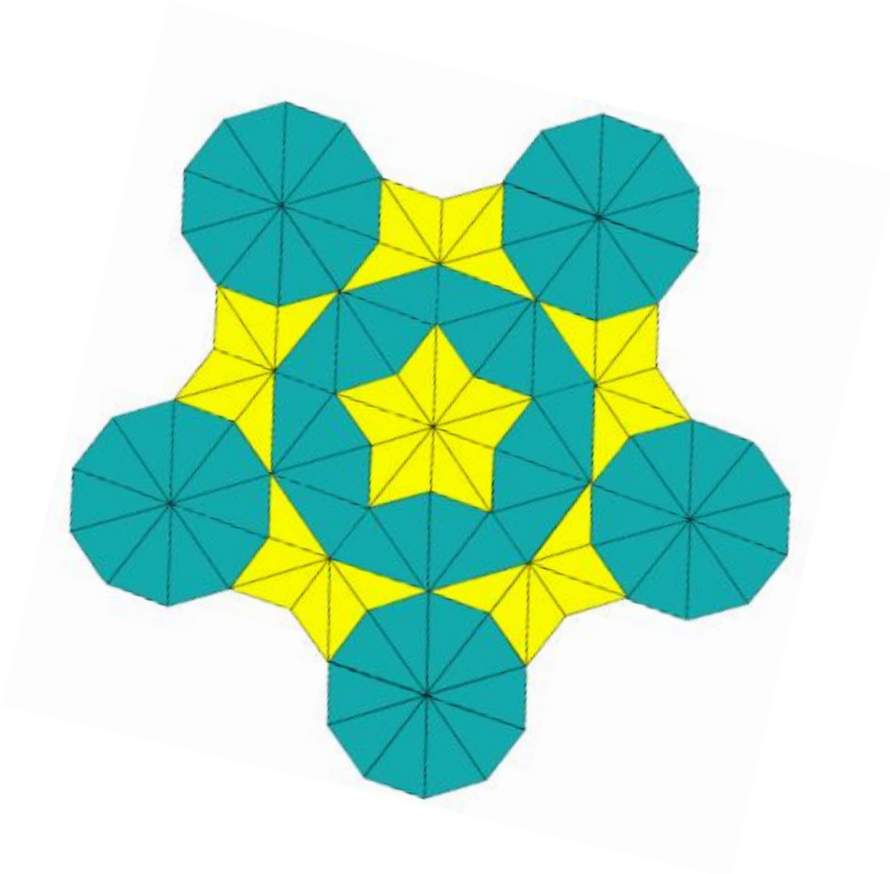

**Figure 9**
A 5-fold symmetric patch of Penrose tiling with darts and kites.

There are two inflation techniques using the Robinson triangles for aperiodic tiling. One version is the Penrose-Robinson tiling (PRT) and the other is the Tubingen triangle tiling (TTT). For further discussions see (Baake & Grimm, 2013).
In appendix A we have displayed the detailed relations between the Voronoi cells of $A_4$ and 5-dimensional cubic lattice.

### 4.3a. Projection of the Voronoi cell of the root lattice $A_5$

Here the Coxeter number is $h = 6$ and the dihedral subgroup $D_6$ of the Coxeter-Weyl group $W(a_5)$ is of order 12. The Delone polytopes are three different orbits $(\omega_1)_{a_5} = -(\omega_5)_{a_5}, (\omega_2)_{a_5} = -(\omega_4)_{a_5}, (\omega_3)_{a_5}$. The Voronoi cell $V(0)$ of the lattice $A_5$ is a polytope in 5-dimensional space and is the disjoint union of the above Delone polytopes. Its 4-dimensional facets are the 4-dimensional rhombohedra implying that the 2-dimensional facets are the rhombuses generated by the pair of vectors $(k_i,\ k_j) = -\frac{1}{6},\ i \neq j = 1, 2, \ldots, 6$. In the Coxeter plane the scalar product would read

$$\left((k_i)_{\parallel}, (k_j)_{\parallel}\right) = \cos\left(\frac{2\pi(j-i)}{6}\right). \tag{21}$$

All possible projected rhombuses turn out to be the one with interior angles $(60^0, 120^0)$. The tiling with this rhombus is known in the literature as the *rhombille* tiling which is used for the study of spin structures in diatomic molecules based on the Ising models. The rhombille tiling is depicted in Fig. 10.

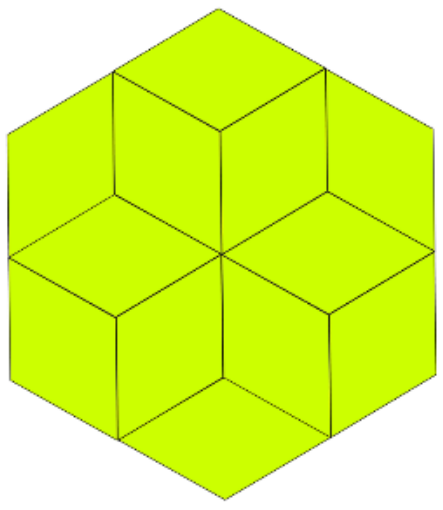

**Figure 10**
The rhombille tiling with a $D_6$ symmetry.

**4.3b. Projection of the Delone cells of the root lattice $A_5$**

The prototiles obtained from the 2-dimensional Delone faces are of three types of triangles with angles $\left(\frac{\pi}{3}, \frac{\pi}{3}, \frac{\pi}{3}\right)$, $\left(\frac{\pi}{6}, \frac{\pi}{6}, \frac{2\pi}{3}\right)$, $\left(\frac{\pi}{6}, \frac{\pi}{3}, \frac{\pi}{2}\right)$ represented by dark blue, yellow and green triangles respectively. Some patches of 6-fold symmetric tilings by three triangles are shown in Fig. 11 which follow from substitution technique.

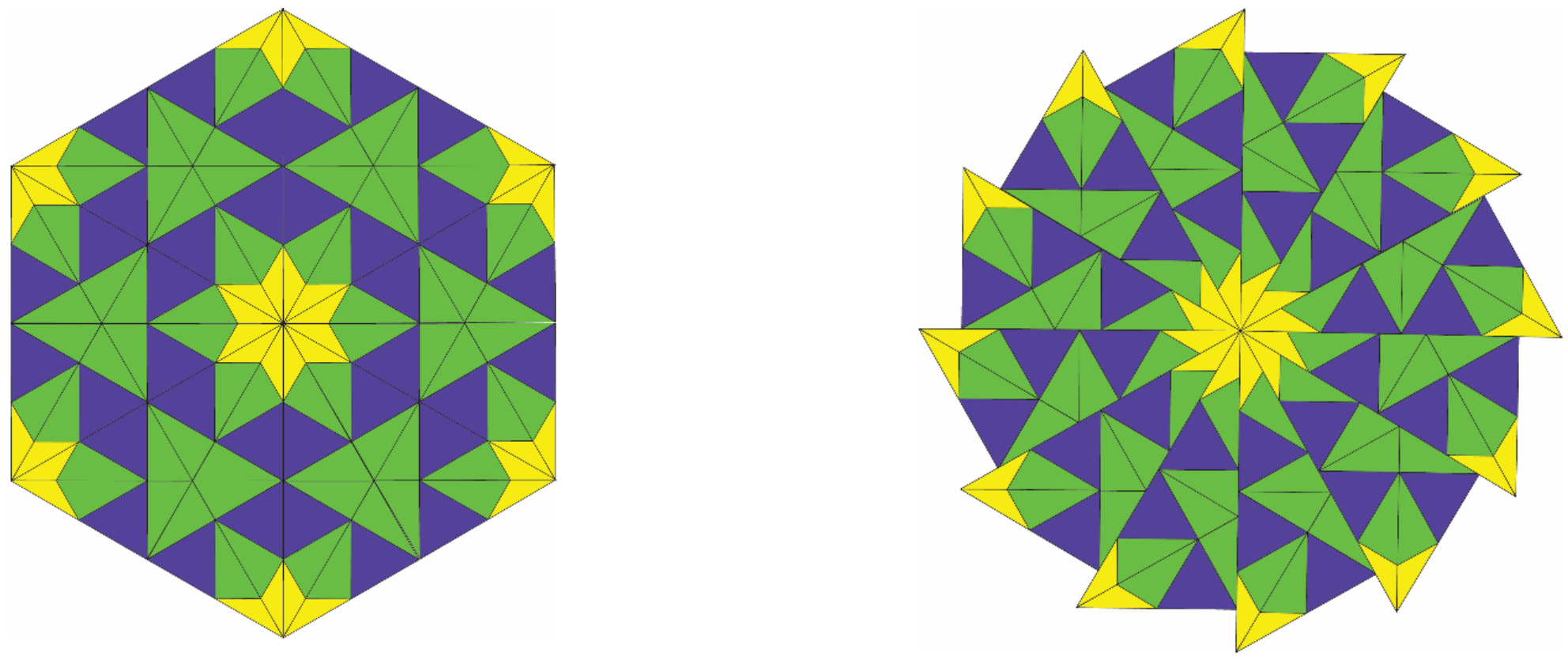

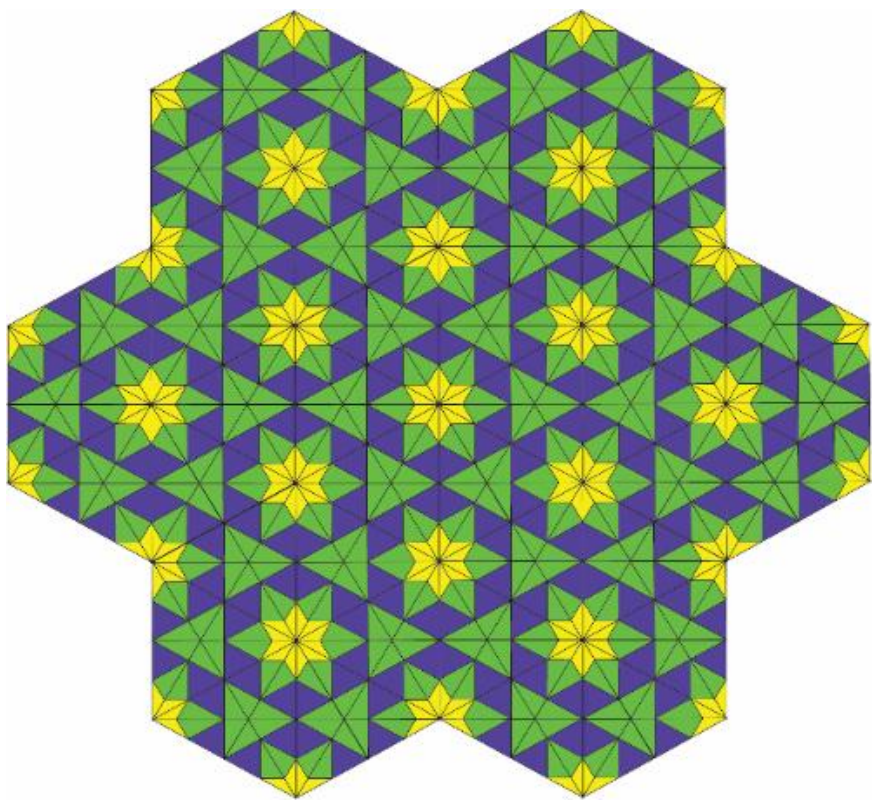

**Figure 11**
Some patches from three triangular prototiles from Delone cells of the root lattice $A_5$ obtained by substitution rule.

### 4.4a. Prototiles from the projection of the Voronoi cell $V(0)$ of the root lattice $A_6$

The rhombic prototiles from the Voronoi cell $V(0)$ can be obtained from the formula $((k_i)_{\|}, (k_j)_{\|}) = \cos\left(\frac{2\pi(j-i)}{7}\right), i \neq j = 1, 2, \ldots, 7$. This would lead to three types of rhombi with angles $\left(\frac{\pi}{7}, \frac{6\pi}{7}\right), \left(\frac{2\pi}{7}, \frac{5\pi}{7}\right), \left(\frac{3\pi}{7}, \frac{4\pi}{7}\right)$. There are numerous discussions on the 7- fold symmetric rhombic tiling with three rhombi. One particular construction is based on the cut-and-projection technique from 7-dimensional cubic lattice (Whittaker& Whittaker, 1988). See also the website of Professor Gerard t' Hooft. A patch obtained by substitution technique is illustrated in Fig. 12.

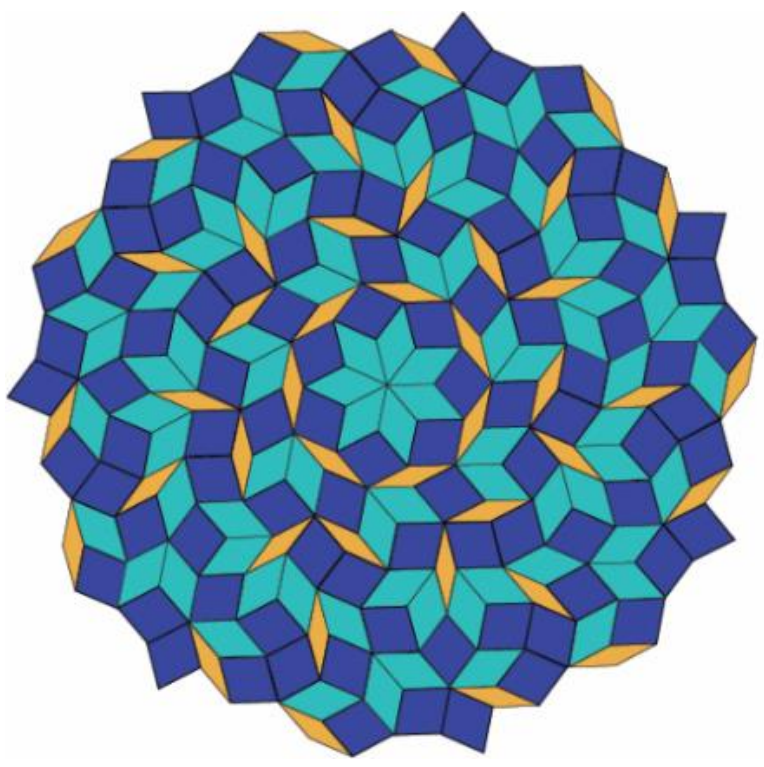

**Figure 12**
A patch of 7-fold symmetric rhombic tiling.

### 4.4b. Projection of the Delone cells of the root lattice $A_6$

The prototiles obtained from 2-dimensional Delone faces are of four types of triangles with angles $\left(\frac{\pi}{7}, \frac{\pi}{7}, \frac{5\pi}{7}\right), \left(\frac{2\pi}{7}, \frac{2\pi}{7}, \frac{3\pi}{7}\right), \left(\frac{3\pi}{7}, \frac{3\pi}{7}, \frac{\pi}{7}\right)$ and $\left(\frac{\pi}{7}, \frac{2\pi}{7}, \frac{4\pi}{7}\right)$. The substitution tilings by the last three triangles are widely known as Danzer tiling. The number of prototiles obtained by projection here is the same as the tiles obtained from the tangents of the deltoid (Nischke & Danzer, 1996). They have illustrated the patches involving only three of the triangles. In Fig.13 we illustrate a patch of 7-fold symmetric tiling comprising of all four triangles.

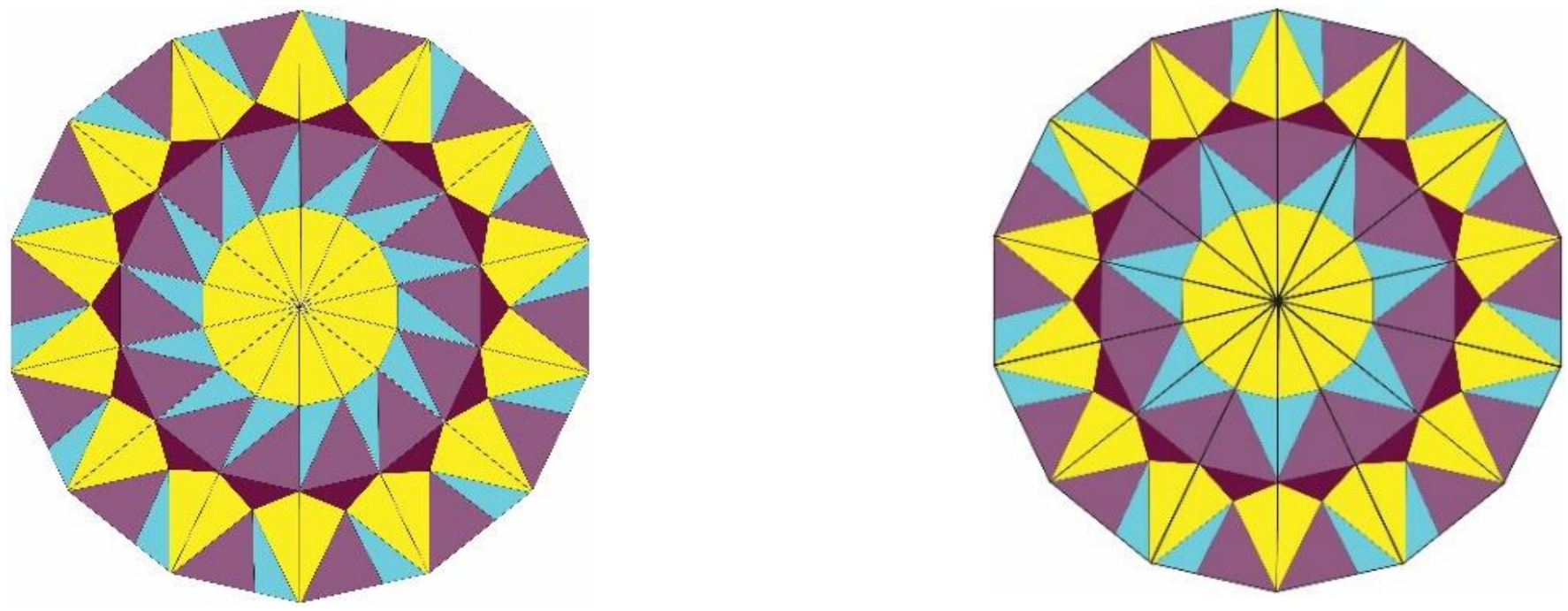

**Figure 13**
Two 7-fold symmetric patches obtained from the projection of the Delone cells of the root lattice $A_6$ .

### 4.5a. Prototiles from the projection of the Voronoi cell $V(0)$ of the root lattice $A_7$

The famous 8-fold symmetric rhombic tiling by two prototiles is known as Amman-Beenker tiling (Grünbaum & Shephard, 1987). Interestingly enough we obtain the same prototiles from the projection of the Voronoi cell $V(0)$ of the root lattice $A_7$ . The point group of the aperiodic tiling is the dihedral group $D_8$ of order 16 and the prototiles are generated by the pair of vectors satisfying the scalar product $\left((k_i)_{\parallel},(k_j)_{\parallel}\right) = \cos\left(\frac{2\pi(j-i)}{8}\right), i \neq j = 1, 2, \ldots, 8$. We obtain two Amman-Beenker prototiles one rhombus with angles $\left(\frac{\pi}{4},\frac{3\pi}{4}\right)$ and a square $\left(\frac{\pi}{2},\frac{\pi}{2}\right)$. A patch of aperiodic tiling is illustrated in Fig. 14.

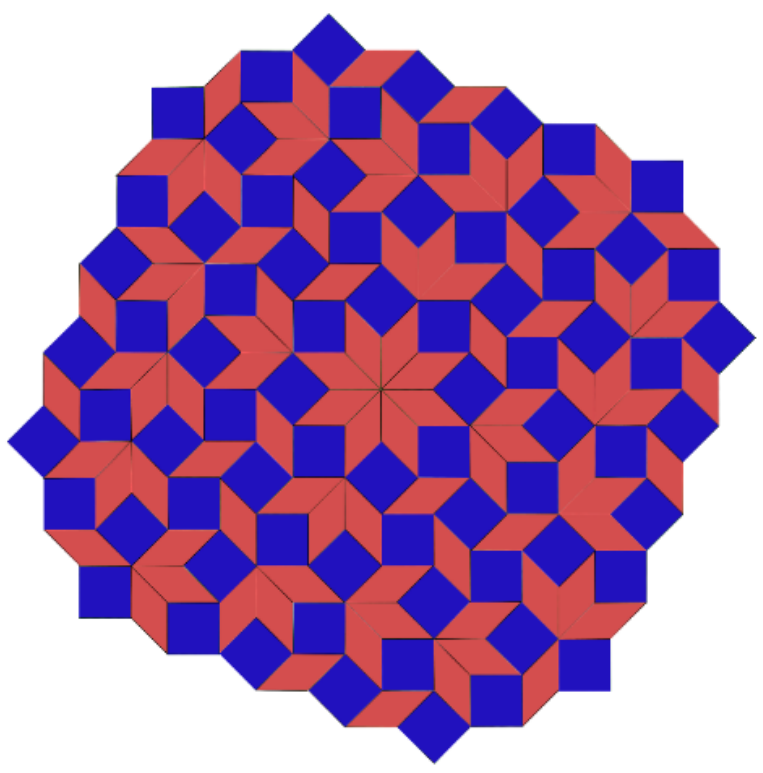

**Figure 14**
The patch obtained from the projection of the Voronoi cell of the root lattice $A_7$ .

### 4.5b. Projection of the Delone cells of the root lattice $A_7$

Substituting $n = 7$ in (13) we obtain five different triangular prototiles from the projection of the Delone cells of the root lattice $A_7$ . Three of the prototiles are the isosceles triangles with angles $\left(\frac{\pi}{8},\frac{\pi}{8},\frac{6\pi}{8}\right)$, $\left(\frac{2\pi}{8},\frac{2\pi}{8},\frac{4\pi}{8}\right)$, $\left(\frac{3\pi}{8},\frac{3\pi}{8},\frac{2\pi}{8}\right)$ and the other two are the scalene triangles with angles $\left(\frac{\pi}{8},\frac{2\pi}{8},\frac{5\pi}{8}\right)$, $\left(\frac{\pi}{8},\frac{3\pi}{8},\frac{4\pi}{8}\right)$. So far as we know no one has studied the tilings of the plane with these prototiles. Two patches of 8-fold symmetric tilings with these prototiles are depicted in Fig.15.

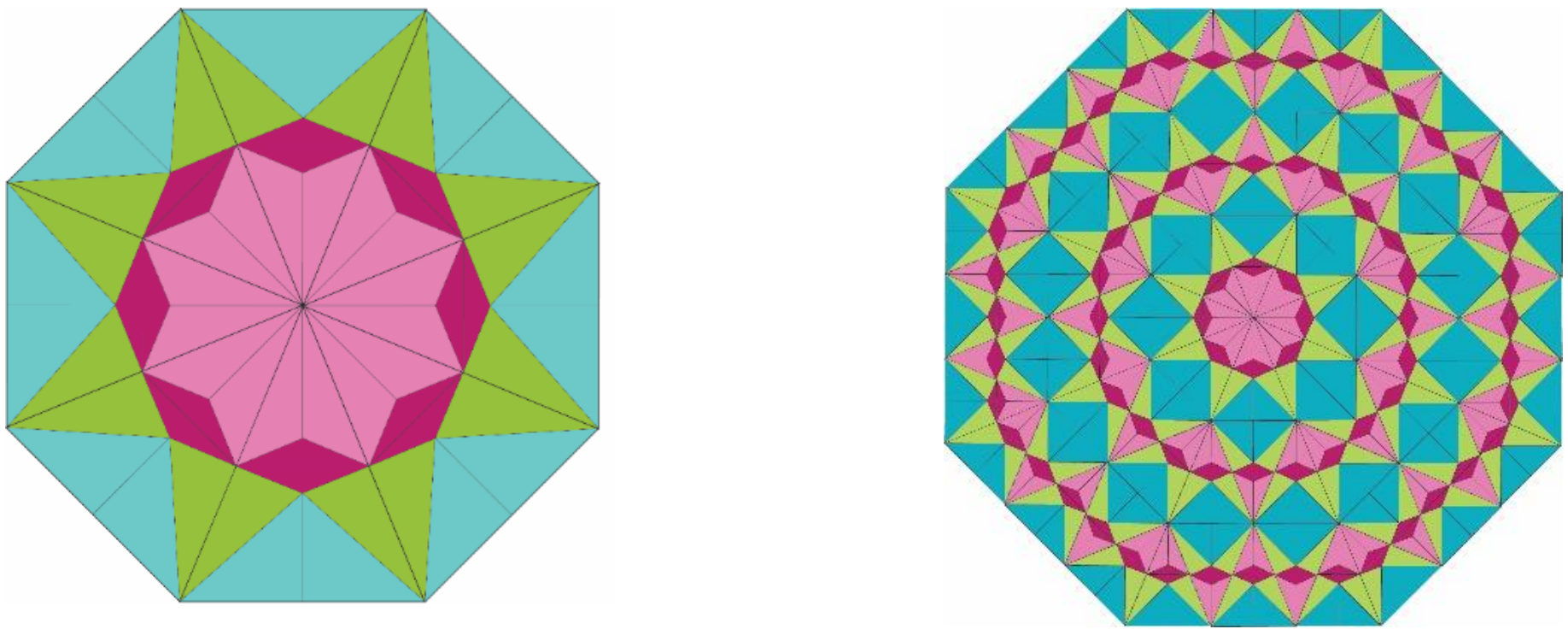

**Figure 15**
Patches of prototiles from Delone cells of the root lattice $A_7$ with 8-fold symmetry.

### 4.6a. Prototiles from projection of the Voronoi cell of the root lattice $A_{11}$

The 12-fold symmetric rhombic tiling has 3 prototiles of rhombuses with interior angles $\left(\frac{\pi}{6},\frac{5\pi}{6}\right)$, $\left(\frac{\pi}{3},\frac{2\pi}{3}\right)$, and $\left(\frac{\pi}{2},\frac{\pi}{2}\right)$. A patch of 12-fold symmetric tiling obtained by substitution technique depicted in Fig.16.

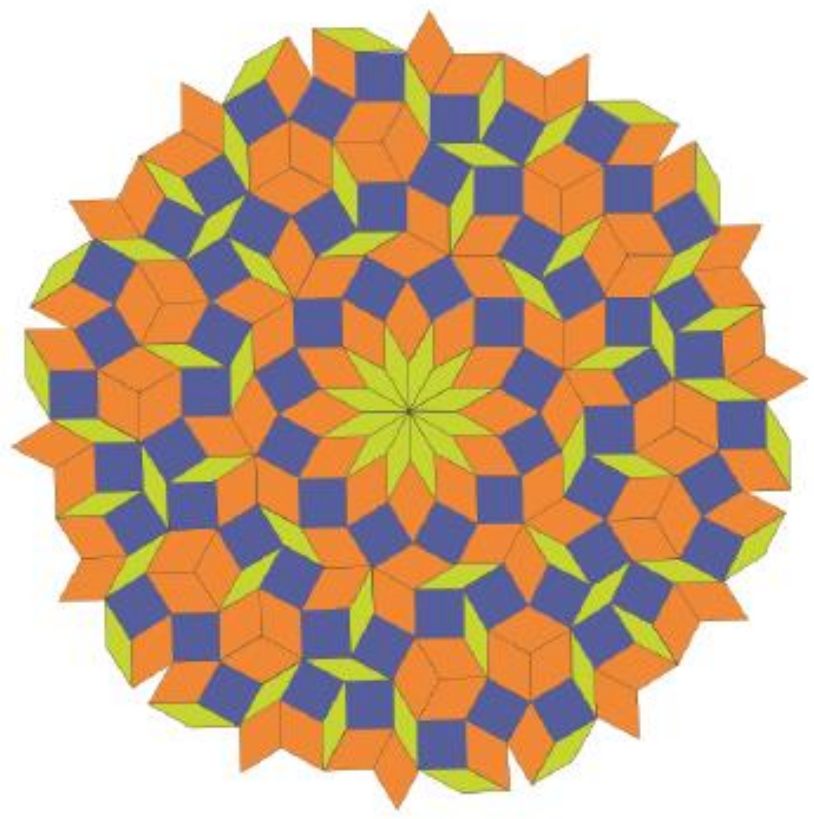

**Figure 16**
Three rhombuses illustrated with different colors obtained from the Voronoi cell of $A_{11}$ and a patch is depicted.

### 4.6b. Prototiles from the projection of the Delone cells of the root lattice $A_{11}$

The number of prototiles in this case is 12. As the rank of the Coxeter-Weyl group is increasing, the number of the triangular prototiles are also increasing. The 12 triangles are given as the set of integers $(n_1, n_2, n_3)$ defined in Sec. 3 as $(1,1,10), (1,2,9), (1,3,8), (1,4,7), (1,5,6), (2,2,8), (2,3,7), (2,4,6), (2,5,5),$ $(3,3,6), (3,4,7), (4,4,4)$. A patch of 12-fold symmetric tiling is depicted in Fig. 17.

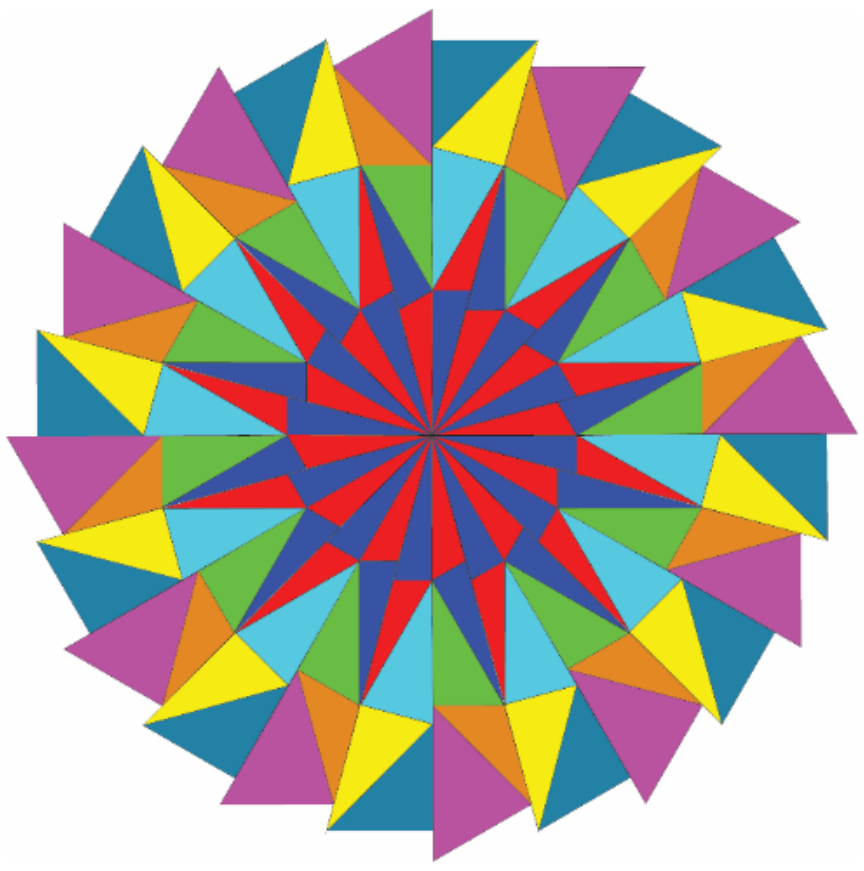

**Figure 17**
A 12-fold symmetric patch obtained by substitution from triangular prototiles originating from Delone cells of the root lattice $A_{11}$.

## 5. Concluding Remarks

Rhombic prototiles usually arise from projection of the higher dimensional cubic lattices $B_{n+1}$ because 2-dimensional square faces project onto rhombuses and represent the local $2(n+1)$-fold symmetric rhombic aperiodic tilings of the Coxeter plane. Depending on the shift of the Coxeter plane one may reduce the tilings to $(n+1)$-fold symmetric aperiodic tilings. The rhombic $(n+1)$-fold symmetric aperiodic tilings can also be obtained from the lattice $A_n$ for the latter lattice is a sub-lattice of the lattice $B_{n+1}$. One may visualize that the lattice $B_{n+1}$ projects onto the lattice $A_n$ at first stage as the Voronoi cell of the cubic lattice projects onto the Voronoi cell of the root lattice $A_n$. So, the classification of the rhombic aperiodic tilings of both lattices are the same. An advantage of the root lattice $A_n$ is that it can be tiled by the Delone cells with triangular 2-dimensional faces. Projections of the Delone cells of the root lattice $A_n$ allow the classifications of the triangular aperiodic tilings. A systematic study of the projections of the Delone and Voronoi cells of the root and weight lattices of the simply laced ADE Lie algebras may lead to more interesting prototiles and aperiodic tilings. The present paper was concerned only with the aperiodic rhombic and triangular tilings of the root lattice exemplifying the 5-fold, 8-fold and 12-fold symmetric aperiodic tilings as they represent some quasicrystallographic structures.

## Appendix A:

### 3-dimensional cube and the Voronoi cell $V(0)$ of $A_2$

Vertices of the Voronoi cell of a 3-dimensional cubic lattice are given by the vectors $\frac{1}{2}(\pm l_1 \pm l_2 \pm l_3)$. Since the vector representing the space diagonal projects to the origin $\pm\frac{1}{2}l_0 \rightarrow 0$ the images of the 6 remaining vertices of the cube can be represented by the vectors $\pm k_i$, $(i = 1, 2, 3)$ in $2D$ space. Each set of vectors $k_i$ and $-k_i$ represents an equilateral triangle. The union of two sets is the hexagon representing the Voronoi cell $V(0)$ of the root lattice $A_2$. In particle physics, it is used to define the color and anti-color degrees of freedom of the quarks. Vertices of the six

Delone cells (equilateral triangles) centralizing the vertices of the Voronoi cell $V(0)$ consists of the vectors $0, \pm(k_i - k_j)$.

**5-dimensional cube and the Voronoi cell $V(0)$ of $A_4$**

Vertices of the Voronoi cell of a 5-dimensional cubic lattice are given by the vectors $\frac{1}{2}(\pm l_1 \pm l_2 \pm l_3 \pm l_4 \pm l_5)$. The cube is cut by hyperplanes orthogonal to the vector $l_0$ at four levels which splits the cube into four Delone polytopes of $A_4$, the union of which represents the Voronoi cell of $A_4$. When projected into 4-dimensional space they can be represented as $(\omega_1)_{a_4} = -(\omega_4)_{a_4}, (\omega_2)_{a_4} = -(\omega_3)_{a_4}$. Vertices of a typical 2-dimensional face of the 5-dimensional cube can be taken as $\frac{1}{2}(l_1 + l_2 + l_3 + l_4 + l_5)$, $\frac{1}{2}(-l_1 + l_2 + l_3 + l_4 + l_5)$, $\frac{1}{2}(l_1 - l_2 + l_3 + l_4 + l_5)$, $\frac{1}{2}(-l_1 - l_2 + l_3 + l_4 + l_5)$. The edges of the square are represented by the vectors parallel to $l_1$ and $l_2$. In 4-dimensional space $l_0$ projects into the origin by (6) and therefore $l_1$ and $l_2$ are replaced by the vectors $k_1$ and $k_2$ constituting a rhombus with obtuse angle $\theta = -\cos^{-1}\left(\frac{1}{4}\right) \cong 104.5^0$. Each 3-dimensional cubic facet of 5-dimensional cube projects into a rhombohedral facet of the Voronoi cell $V(0)$ as it is shown in the following example.

The set of vertices $\frac{1}{2}(l_1 \pm l_2 \pm l_3 \pm l_4 - l_5)$ is a 3-dimensional cube centered at the point $\frac{1}{2}(l_1 - l_5)$. This cube projects into a 3-dimensional rhombohedral facet of the Voronoi cell $V(0)$ of $A_4$ as shown in Fig. 6 where the vertices of the rhombohedron are given by

$$\begin{aligned}
&\frac{1}{2}(l_1 - l_2 - l_3 - l_4 - l_5) \to k_1\\
&\frac{1}{2}(l_1 + l_2 + l_3 + l_4 - l_5) \to -k_5,\\
&\frac{1}{2}(l_1 + l_2 - l_3 - l_4 - l_5) \to k_1 + k_2,\\
&\frac{1}{2}(l_1 - l_2 + l_3 - l_4 - l_5) \to k_1 + k_3,\\
&\frac{1}{2}(l_1 - l_2 - l_3 + l_4 - l_5) \to k_1 + k_4,\\
&\frac{1}{2}(l_1 + l_2 + l_3 - l_4 - l_5) \to -(k_4 + k_5),\\
&\frac{1}{2}(l_1 + l_2 - l_3 + l_4 - l_5) \to -(k_3 + k_5),\\
&\frac{1}{2}(l_1 - l_2 + l_3 + l_4 - l_5) \to -(k_2 + k_5).
\end{aligned} \tag{A.1}$$

All the rhombohedral facets of the Voronoi cell $V(0)$ are obtained by permutations. Centers of the rhombohedral facets are the halves of the roots of the root system. This also explains why the Voronoi cell of $A_4$ is the dual polytope of the root polytope.

In brief, the square face of 5-dimensional cube projects onto rhombic face of the Voronoi cell $V(0)$ in 4-dimensional space and all 3-dimensional cubic facets project onto the 3-dimensional rhombohedra. Of course, the result of both projections leads to the same Penrose rhombic aperiodic tilings. This is true for all $B_{n+1} \rightarrow A_n$ projections, the only difference is that the aperiodic tiles from the projections of $B_{n+1}$ exhibit $2(n+1)$-fold symmetry while $A_n$ projections display $(n+1)$-fold symmetry only.

**References**

Baake, M., Joseph, D., Kramer, P. & Schlottmann, M. (1990). *J. Phys. A: Math. & Gen.* **23**, L1037-L1041.
Baake, M. & Grimm, U. (2013). *Aperiodic Order, Volume 1: A Mathematical Invitation*, Cambridge University Press, Cambridge.
Boyle, L. & Steinhardt, P. J. (2016). *ArXiv*: 1608.08215.
Carter, R. W. (1972). *Simple Groups of Lie Type*, pp:158-169. John Wiley & Sons, New York.
Chen, L., Moody, R. V. & Patera, J. (1998). *Fields Institute for Research in Mathematical Sciences Monographs Series* 10, pp. 135-178. Providence, Rhode Island, USA: AMS.
Conway, J. H. & Sloane, N. J. A. (1988). *Sphere Packings, Lattices and Groups*. Springer-Verlag, New York Inc.
Conway, J. H. & Sloane, N. J. A. (1991). *Miscellanea Mathematica*, (ed. Hilton, P., Hirzebruch, F. & Remmert, R.), pp. 71-108. New York: Springer.
Coxeter, H. S. M. (1973). *Regular Polytopes*, Third Edition, Dover Publications.
de Bruijn, N. G. (1981). *Nederl. Akad. Wetensch. Proceedings Ser.***A84** (=Indagationes Math. **43**), 38-66.
Delaunay, N. B. (1929). *Izv. Akad. Nauk SSSR Otdel. Fiz.-Mat. Nauk*. pp. 79-110 and 145-164.
Delaunay, N. B. (1938a). *Usp. Mat. Nauk*. **3,** 16-62.
Delaunay, N. B. (1938b). *Usp. Mat. Nauk*. **4,** 102-164.
Deza, M. & Grishukhin, V. (2004). *Geometriae Dedicate* **104**, 15-24.
Di Vincenzo, D. & Steinhardt, P. J. (1991). *Quasicrystals: the state of the art*, World Scientific Publishers, Singapore.
Duneau, M. & Katz, A. (1985). *Phys. Rev. Lett.* **54**, 2688–2691.
Engel, P. (1986). *Geometric Crystallography: An Axiomatic Introduction to Crystallography*. Dordrecht: Springer.
Engel, P., Michel, L. & Senechal, M. (1994). *preprint* IHES /P/04/45.
Grünbaum, B. (1967). *Convex Polytopes*, Wiley, New York.
Grünbaum, B. & Shephard, G. C. (1987). *Tilings and Patterns*, Freeman, New York.
Humphreys, J. E. (1990). *Cambridge Studies in Advanced Mathematics, Reflection Groups and Coxeter Groups*, Vol. 29. Cambridge University Press.
Janot, C. (1993). *Quasicrystals: a primer*. Oxford University Press.
Koca, N. O., Koca, M. & Al-Siyabi, A. (2018a). *Int. J. Geom. Methods Mod. Phys.***15** (4), 1850058.
Koca, M., Ozdes Koca, N., Al-Siyabi, A. & Koc, R. (2108b). *Acta Cryst*. A**74**, 499-511.
Koca, M., Koca, N. O. & Koc, R. (2014). *Int. J. Geom. Meth. Mod. Phys*. **11**, 1450031
Koca, M., Koca, N. & Koc, R. (2015). *Acta Cryst* A**71**, 175-185.
Lagarias, J. C. (1996). *Commun. Math. Phys*. **179**, 365-376.
Meyer, Y. (1972). *Algebraic Numbers and Harmonic Analysis*, North Holland, Amsterdam.

Michel, L. (1995). *Bravais classes, Voronoi cells, Delone cells, Symmetry and Structural Properties of Condensed Matter*, (ed. Lulek, T., Florek, W. & Walcerz, S., Zajaczkowo), World Sci., Singapore, p. 279.
Michel, L. (2001). *Complete description of the Voronoi cell of the Lie algebra $A_n$ weight lattice. On the bounds for the number of d-faces of the n-dimensional Voronoi cells, Algebraic Methods in Physics*, Yvan Saint-Aubin and Luc Vinet (editors), Springer, pp. 149-171.
Moody, R.V. & Patera, J. (1992). *J. Phys. A: Math. Gen.* **25**, 5089-5134.
Moody, R. V. (1997). *Meyer Sets and their duals, in The Mathematics of Long-Range Aperiodic Order*, *NATO ASI Series C* **489**, pp 403-441. Kluwer Academic Publishers, The Netherlands.
Nischke, K. P. & Danzer, L. (1996). *Discrete Comput. Geom.* **15**, 221-236.
Penrose, R. (1974). *Bull. Inst. Math. Appl*. **10**, 266-271.
Penrose, R. (1978). *Eureka* **39**, 16-22.
Shechtman, D., Blech, I., Gratias, D. & Cahn, J.W. (1984). *Phys. Rev. Lett*. **53**, 1951-1953.
Senechal, M. (1995). *Quasicrystals and Geometry*, Cambridge University Press, Cambridge.
Voronoi, G. (1908), *Journal für die Reine und Angewandte Mathematik*, **134**, 198-287
Voronoi, G. (1909). *Journal für die Reine und Angewandte Mathematik*, **136**, 67-181.
Whittaker, E. J. W. & Whittaker, R. M. (1987). *Acta Cryst*. **A44**, 105-112.